\documentclass{article}

\usepackage{amssymb}
\usepackage{amsmath}
\usepackage{graphics}

\let\<\langle
\let\>\rangle
\font\bigcmsy=cmsy10.pk scaled 2000
\def\bigtimes{\mathop{\,\vrule width0pt depth2pt height8pt
            \smash{\lower2pt\hbox{\bigcmsy\char'002}}\,}\limits}

\begin{document}

\begin{center}
\Large{\textbf{Generalization of order separability for free
products and omnipotence of free products of groups.}}
\end{center}
\begin{center}

\textbf{Vladimir V. Yedynak}

\textsl{Faculty of Mechanics and Mathematics, Moscow State
University}

\textsl{Moscow 119992, Leninskie gory, MSU}

\textsl{edynak\_vova@mail.ru}
\end{center}

\begin{abstract}

It was proved that for any finite set of elements of a free product
of residually finite groups such that no two of them belong to
conjugate cyclic subgroups and each of them do not belong to a
subgroup which is conjugate to a free factor there exists a
homomorphism of the free product onto a finite group such that the
order of the image of each fixed element is an arbitrary multiple of
a constant number.


\end{abstract}

\section{Introduction}

Order separabilities are connected with the investigation of the
correlation between the orders of elements' images after a
homomorphism of a group onto a finite group. For example in [1] it
was proved that for each elements $u$ and $v$ of a free group $F$
such that $u$ is conjugate to neither $v$ nor $v^{-1}$ there exists
a homomorphism of $F$ onto a finite group such that the images of
$u$ and $v$ have different orders. In [6] it was proved that this
property is inherited by free products. This paper is devoted to the
proof of the theorem that strengthens the property of order
separability for the class of free products of groups.

\textbf{Theorem.} Consider the group $G=A\ast B$ where the subgroups
$A$ and $B$ are residually finite. Consider the elements
$u_1,\ldots, u_n$ such that $u_i\in G\setminus\{\,\cup_{g\in
G}(g^{-1}Ag\cup g^{-1}Bg)\}, u_i, u_j$ belongs to conjugate cyclic
subgroups whenever $i=j$. Then there exists the natural number $K$
such that for each ordered sequence $l_1,\ldots, l_n$ of natural
numbers there exists a homomorphism $\varphi$ of $G$ onto a finite
group such that the order of $\varphi(u_i)$ is equal to $Kl_i$

The property under study in this work is closely connected with
omnipotence which was investigated in [2], [3] where it was shown
that free groups and fundamental groups of compact hyperbolic
surfaces are omnipotent. Besides all finite sets of independent
elements whose orders are infinite in a Fuchsian group of the first
type also satisfy the property of omnipotence [4].

Definition. The group $G$ is called omnipotent if for each elements
$u_1,\ldots, u_n$ such that no two of them have conjugate nontrivial
powers there exists a number $K$ such that for each ordered sequence
of natural numbers $l_1,\ldots, l_n$ there exists a homomorphism
$\varphi$ of $G$ onto a finite group such that the order of
$\varphi(u_i)$ equals $Kl_i$.

The familiar property was also investigated in [7] where some
sufficient conditions were found for $n$-order separability of free
products. The group $G$ is said to be $n$-order separable if for a
set $S=\{\,s_1,\ldots, s_n\mid s_i\neq h^{-1}s_j^{\pm1}h, i\neq j\}$
of $n$ elements of $G$ there exists a homomorphism of $G$ onto a
finite group mapping $S$ onto a set whose elements have pairwise
different orders.

Notice that the theorem of this paper will enable to investigate the
residual properties of the fundamental group of graphs of groups
whose vertex groups are residually finite free products and edge
groups are cyclic not belonging to subgroups conjugate to free
factors of vertex groups.

\section{Notations and Definitions}

We consider that for every graph there exists a mapping $\eta$ from
the set of edges of this graph onto itself. For every edge $e$ this
mapping corresponds an edge which is inverse to $e$. Besides the
following conditions are true: $\eta(\eta(e))=e$ for each $e$,
$\eta$ is a bijection, for every edge $e$ the beginning of $e$
coincides with the end of the edge $\eta(e)$.

The graph is called oriented if from every pair of mutually inverse
edges one of them is fixed. The fixed edge is called positively
oriented and the inverse edge is called negatively oriented.

Let $G$ be a free product of groups $A$ and $B$. There exists a
correspondence such that for every action of $G$ on the set $X$ at
which both $A$ and $B$ act freely there exists a graph $\Gamma$
satisfying the following properties:

1) for each $c\in A\cup B$ and for each vertex $p$ of $\Gamma$ there
exists exactly one edge labelled by $c$ going into $p$ and there
exists exactly one edge labelled by the element $c$ which goes away
from $p$.

2) for every vertex $p$ of $\Gamma$ the maximal connected subgraph
$A(p)$ of $\Gamma$ containing $p$ whose positively oriented edges
are laballed by the elements of $A$ is the Cayley graph of the group
$A$ with generators $\{\,A\}$; we define analogically the subgraph
$B(p)$.

3) we consider that for every edge $e$ from the first item there
exists the edge inverse to $e$ which does not bear a label; two
edges with labels are not mutually inverse; edges with labels are
positively oriented.

Definition 1. We say that a graph is the free action graph of the
group $G=A\ast B$ if it satisfies the properties 1), 2), 3).

Note that if $\varphi$ is the homomorphism of the group $G$ such
that $\varphi_{A\cup B}$ is the bijection then the Cayley graph
$Cay(\varphi(G); \{\,\varphi(A)\cup\varphi(B)\})$ of the group
$\varphi(G)$ with respect to the set of generators
$\{\,\varphi(A)\cup\varphi(B)\}$ is the free action graph of the
group $G$.

Remark. In what follows appending a new edge with label to a free
action graph we shall consider that it is positively oriented and
the inverse edge would have been appended. And if we delete an edge
with label the inverse edge would have been deleted.

If $e$ is the edge then $\alpha(e), \omega(e)$ are vertices which
coincide with the beginning and the end of $e$ correspondingly.

If we have the free action graph $\Gamma$ of the group $G$ then
there exists the action of $G$ on the set of vertices of $\Gamma$
which is defined as follows. Let $p$ be an arbitrary vertex of
$\Gamma$. Then according to the definition of the free action graph
for each element $c$ from $A\cup B$ there exist edges $e$ and $f$
whose labels are equal to $c$ such that $\alpha(e)=p, \omega(f)=p$.
In this case the action of $c$ on $p$ is defined as follows: $p\circ
c = \omega(e), p\circ c^{-1} = \alpha(f)$.

Remark also that if we change the property 2) in the definition of
the free action graph supposing that $A(p)$ and $B(p)$ are the
Cayley graphs of the homomorphic images of the groups $A$ and $B$
correspondingly we also obtain the graph such that there exists the
action of the group $G$ on the set of its vertices. Such a graph
will be referred to as an action graph of the group $G$.

Since there exists the action of $G$ on the set of vertices of an
action graph $\Gamma$ there exists a homomorphism of $G$ onto the
group $S_n$, where $n$ is the cardinal number of the set of vertices
of the graph $\Gamma$. Having a group $G$ and its action graph
$\Gamma$ we shall denote this homomorphism as $\varphi_{\Gamma}$.

If $e$ is the positively oriented edge of the action graph, then
Lab$(e)$ is the label of $e$.

Definition 2. Let $u$ be a cyclically reduced element of the group
$G$ which belongs to neither $A$ nor $B$ and $\Gamma$ is the action
graph of $G$. Fix a vertex $p$ of $\Gamma$. Then $u$-cycle in this
action graph going from $p$ is the cycle $R=e_1\ldots e_n$ which
satisfies the following properties:

1) the path $P$ is a closed path such that its beginning
$\alpha(P)=p$

2) consider $u=u_1\ldots u_k$ where $u_i\in A\cup B, u_i, u_{i+1}$
as well as $u_1, u_k$ do not belong to one free factor
simultaneously; then $k$ divides $n$ and the edge $e_{ik+j}$ is
positively oriented and has a label $u_j$, $1\leqslant j\leqslant k$
(indices are modulo $n$)

3) the cycle $P$ is the minimal cycle which satisfies properties 1),
2).

Definition 3. Suppose we have a path $S=e_1\cdots e_n$ in the action
graph. Then the label of this path is the element of the group which
is equal to $\prod_{i=1}^n$Lab$(e_i)'$, where Lab$(e_i)'$ equals
either the label of $e_i$, if this edge is positively oriented, or
Lab$(e_i)'= $ Lab$(\eta(e_i))^{-1}$ otherwise. We shall denote the
label of the path $S$ as Lab$(S)$.

Definition 4. Fix the graph $\Gamma$, $p$ and $q$ are vertices from
$\Gamma$. Then we define the distance between $p$ and $q$ as
$\rho(p, q)=\min_{S}\ l(S)$, where $S$ is an arbitrary path
connecting $p$ and $q$, $l(S)$ is the number of edges in $S$.

Notice that if a cycle $S$ does not have $l$-near vertices then each
subpath of $S$ of length which less or equal than $l$ is geodesic.

Definition 5. Fix an arbitrary graph and a cycle $S=e_1\cdots e_n$
in it. For every nonnegative integer number $l$ we shall say that
$S$ does not have $l$-near vertices, if for every $i, j, i\neq j,
1\leqslant i, j\leqslant n$ the distance between the vertices
$\alpha(e_i), \alpha(e_j)$ is greater or equal than $\min(l+1,
|i-j|, n-|i-j|)$.

Definition 6. Suppose we have the $u$-cycle $S$. It is obvious that
its label equals the $k$-th power of $u$ for some $k$. Then we say
that the length of the $u$-cycle $S$ is equal to $k$.

Note that for the action graph $\Gamma$ and cyclically reduced
element $u\in G\setminus\{\,A\cup B\}$ the order of
$\mid\varphi_{\Gamma}(u)\mid$ coincides with the less common
multiple of lengths of all $u$-cycles in the graph $\Gamma$. Hence
if there exists a $u$-cycle in the action graph whose length equals
$t$ then $\mid\varphi_{\Gamma}(u)\mid$ is a multiple of $t$.

Suppose $u$ is an element of $A\ast B$ and $u=u_1\cdots u_n$ is the
irreducible form of $u$. Then the length of $u$ is the number
$l(u)=n$. The cyclic length $l'(u)$ of an element $u$ is the length
of the cyclically reduced element which is conjugate to $u$.

\section{Auxiliary lemmas}

\textbf{Lemma 1.} Consider the group $G=A\ast B$ where $A$ and $B$
are finite, $l$ and $n$ are natural numbers, $Q$ is a finite set of
elements from $G$ which are cyclically reduced and whose lengths are
greater than 1. Then for each $v\in Q$ there exists the homomorphism
$\varphi$ of $F$ onto a finite group such that for each $u$ from $Q$
the $u$-cycles in the Cayley graph $Cay(\varphi(G);
\{\,\varphi(A)\cup\varphi(B)\})$ of the group $\varphi(G)$ do not
have $l$-near vertices, $\mid\varphi(v)\mid>n$, and $\varphi_{A\cup
B\cup Q}$ --- injection.

\textsl{Proof.}

For each $q\in Q$ define the set $L_q$ which consists of the
elements from $G$ whose length is less or equal than $l+2l(q)+10$
and which do not belong to the subgroup generated by the element
$q$. It is well known that free group are subgroup separable [5].
Besides if a group is virtually subgroup separable than it is
subgroup separable (see [2] for example). Considering the above
there exists the homomorphism $\varphi_q$ of the group $G$ onto a
finite group such that $\varphi_q(L_q)\cap\<\varphi_q(q)\>$ is an
empty set. There also exists the homomorphism $\varphi_v'$ of $G$
onto a finite group such that $\varphi_v'(v^i)\neq1$ where
$i=1,\ldots, n$ and $\varphi_v'|_{A\cup B\cup Q}$ is the injection
since virtually free groups are residually finite. The homomorphism
$\varphi: G\rightarrow(\times_{h\in
Q}\varphi_h(G))\times\varphi_v'(G), \varphi: f\mapsto\prod_{h\in
Q}(\varphi_h(f))\varphi_v'(f)$ is as required. Lemma 1 is proved.

The following statement was proved in [2]

\textbf{Lemma 2.} Consider a group $G$ and its elements $g_1,\ldots,
g_n$ possessing the property that for each $j, 1\leqslant j\leqslant
n,$ there exist constants $K_{j, 1},\ldots, K_{j, n}$ such that for
each natural $m$ there exists a homomorphism $\varphi_{j, m}$ of $G$
to a finite group with the condition that $|\varphi_{j,
m}(g_k)|=K_{j, k}$ for all $k\neq j$ and $|\varphi_{j,
m}(g_j)|=mK_{j, j}$. Then there exists the number $K$ such that for
each ordered sequence of natural numbers $l_1,\ldots, l_n$ there
exists the homomorphism $\psi$ of $G$ onto a finite group satisfying
the property that $|\psi(g_i)|=Kl_i$.

\section{Proof of the theorem}

It follows from lemma 2 that that the theorem can be derived from
the following proposition.

\textbf{Proposition.} Let $G=A\ast B$ be a free product of
residually finite groups $A$ and $B$, $u, v_1,\ldots, v_n \in G$.
Elements $u$ and $v_i$ do not belong to conjugate cyclic subgroups.
Besides $u$ does not belong to a subgroup which is conjugate to
either $A$ or $B$. Then there exist natural numbers $L, K_1,\ldots,
K_n$ such that for each natural $i$ there exists a homomorphism
$\varphi$ of $G$ onto a finite group such that
$\mid\varphi(u)\mid=Li, \mid\varphi(v_i)\mid=K_i, 1\leqslant
i\leqslant n$.

\textsl{Proof.}

Since $A, B$ are residually finite we may consider that $A$ and $B$
are finite. Consider also that the elements $u, v_1,\ldots, v_n$ of
$A\ast B$ are cyclically reduced.

Let us to define the following notation. Consider the action graph
$\Gamma$ of the group $K\ast L$. Let $S$ be the subset of $\Gamma$
(e. g. vertex, edge, path, subgraph etc). Then having a set
$\Gamma_1,\ldots, \Gamma_n$ of copies of $\Gamma$ we consider that
$S^i$ denotes the subset of $\Gamma_i$ corresponding to $S$ in
$\Gamma$.

Put $s=\max_{b\in\{\,u, v_1,\ldots, v_n\}}l(b)$, and let $k'$ be an
arbitrary natural number such that $k'l(u)\geqslant10s$. Put
$k=k'l(u)$. Denote by $P$ the set of all nonunit elements whose
length is less or equal than $10k$. For $Q=\{\,u, v_1,\ldots,
v_n\}\cup P$ according to lemma 1 there exists the homomorphism
$\varphi$ of $G$ onto a finite group such that $\varphi_{A\cup B\cup
Q}$ is the injection and for each $s\in S$ which is cyclically
reduced and whose length is greater than 1 each $s$-cycle in the
graph $\Gamma=Cay(\varphi(G); \{\,\varphi(A)\cup\varphi(B)$ has no
$(k+4)$-near vertices and $\mid\varphi(u)\mid>10k$.

Fix a natural number $m>2$ whose value we shall choose later.
Consider $m$ copies of the graph $\Gamma$: $\Gamma_{i}, 1\leqslant
i\leqslant m$. In the graph $\Gamma$ we fix a $u$-cycle $S=e_1\cdots
e_r$. Without loss of generality we consider that Lab$(e_1)\in A$.
Put $p_i=\alpha(e_i)$, Lab $(e_i)=u_i$ (see Figure 1). For each $i,
1\leqslant i\leqslant m,$ we delete edges incident to $p_2^i$ whose
labels belong to $A$ and delete also edges labelled by the elements
of $A$ whose begin or end points are $p_{k+2}^i$. For each $i$ we
shall denote the obtained graph as $\Gamma_i'$.

\begin{figure}
\centering
\includegraphics[15cm,8cm]{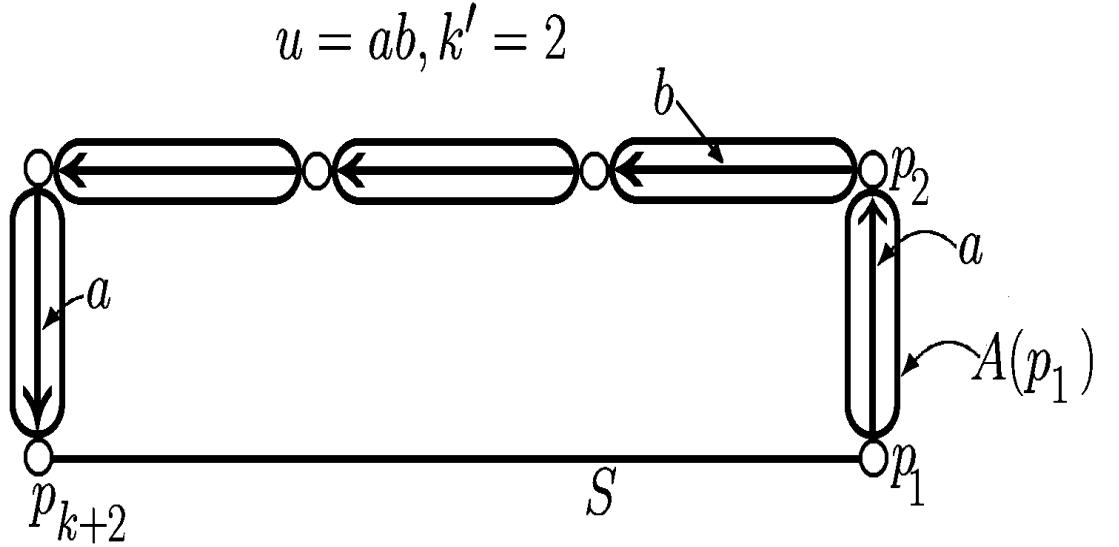}
\caption{The graph $\Gamma$}
\end{figure}

Let $\psi$ be the bijection between the subgraphs $A(p_1)$ and
$A(p_{k+1})$ which saves labels of edges and $\psi(p_1)=p_{k+1}$.

Fix an arbitrary edge $e$ of $\Gamma$ from the subgraph $A(p_1)$
such that the corresponding edge $e^i$ of $\Gamma_i$ was deleted.
Let $q=\alpha(e), r=\omega(e)$. For each $i, 1\leqslant i\leqslant
m,$ if $q\neq p_2$ we connect the vertices $q^i$ and $\psi(r)^{i+1}$
by the new edge $f_i$. If $r\neq p_2$ we connect $r^i$ and
$\psi(q)^{i+1}$ by the edge $f_i$. In both cases the label of $f_i$
coincides with Lab $(e)$, besides if $q\neq p_2$ then $f_i$ goes
away from $q^i$ and if $r\neq p_2$ then $f_i$ goes into $r^i$.

Now we need to complement the structure of obtained graph for to get
the action graph of the group $A\ast B$. But it will not be the free
action graph.

For each $i, 1\leqslant i\leqslant m,$ let us to add one new vertex
$n_i$ to the subgraph $\Gamma_i'$.

Consider an arbitrary edge $e$ from $A(p_{k+2})$ such that the
corresponding edge $e^i$ was deleted from $\Gamma_i$. Put
$q=\alpha(e), r=\omega(e)$. If $q=p_{k+2}$ then connect the vertices
$r^i$ and $n_i$ by the edge $g_i$. If $r=p_{k+2}$ then the new edge
$g_i$ connects the vertices $q^i$ and $n_i$. Put Lab $(g_i)=$ Lab
$(e)$. The begin point of $g_i$ coincides with either $n_i$ or
$q^i$.

For each $c\in A\cup B$ and for each vertex $p$ of the obtained
graph which is not incident to an edge with label $c$ add a loop
with label $c$ going from $p$.

If we fix $j$ then the union of the graph $\Gamma_{j}'$ and
$A(p_{k+1}^j), A(p_{k+2}^j), A(p_1^j)$ is denoted by $\Delta_j$.

We constructed the new graph $\Delta$ which contains subgraphs
$\Gamma_{j}'$ and $\Delta_j$ and is the action graph of the group
$G$.

In the graph $\Delta$ the $u$-cycle $S'$ going from the vertex
$p_{1}^1$ has the length $(\mid\varphi(u)\mid - k')m$. From the
properties of the homomorphism $\varphi$ it follows that
$\mid\varphi(u)\mid>10k=10k'l(u)>k'$. Hence
$\mid\varphi_{\Delta}(u)\mid\geqslant(\mid\varphi(u)\mid - k')m>m$.

Let us to prove that for each $i$ and for each $v_i$-cycle $T$ in
the graph $\Delta$ all vertices of $T$ belong to two subgraphs
$\Delta_{j_1}, \Delta_{j_1+1}$ for some $j_1$. Suppose the contrary.
That is we suppose that there exist pairwise different numbers $j_1,
j_2, j_3$ such that the vertices of $T$ belong to all three
subgraphs $\Delta_{j_1}, \Delta_{j_2}, \Delta_{j_3}$. Note that
different subgraphs $\Delta_{k_1}, \Delta_{k_2}$ has the nonempty
intersection if and only if $\mid k_1-k_2\mid=1$ and their
intersection equals the subgraph $A(p_1^l)$ since $m>2$ where $l$ is
equal to either $k_1$ or $k_2$. So if $\Delta_{j_1}, \Delta_{j_2},
\Delta_{j_3}$ contain vertices of $T$ there exists the number $j$
such that the subgraphs $\Delta_j, \Delta_{j+1}, \Delta_{j+2}$
contain the vertices of $T$ and there exists the path $R$ which is
the part of $T$ and which belongs to
$\Delta_j\cup\Delta_{j+1}\cup\Delta_{j+2}$, $R$ goes away from the
vertex of $\Delta_j$ and goes into the vertex of $\Delta_{j+2}$
(indices are modulo $m$).

From the properties of $R$ it follows that $R$ contains its first
and the last edges $t_j, r_j$ correspondingly such that $t_j\in
A(p_1^j), \omega(t_j)=p_{k+2}^{j+1}, r_j\in A(p_1^{j+1}),
\omega(r_j)=p_{k+2}^{j+2}$, and the rest edges of $R$ are in
$\Delta_{j+1}$.

Because of our supposition that $T$ goes from $\Delta_{j+1}$ into
$\Delta_{j+2}$ it is possible to deduce that $R$ contains the
subpath $s_1\cdots s_l$ such that $s_2\cdots s_{l-1}$ belongs to
$\Gamma_{j+1}'$ and edges $s_1, s_l$ satisfy the following
properties: $\alpha(s_1)\in A(n_{j+1})\cup B(p_{k+2}^{j+1}),
\omega(s_l)\in A(p_1^{j+1})\cup B(p_2^{j+1})$.

Denote the path $e^{j+1}_1e^{j+1}_2\cdots e^{j+1}_{k+1}$ as $S_u$
and $s_2\cdots s_{l-1}$ as $S_{v_i}$ (see Figure 2). Note that
$\rho(\alpha(S_u), \alpha(S_{v_i}))\leqslant1, \rho(\omega(S_u),
\omega(S_{v_i}))\leqslant2$ (the function $\rho$ is taken with
respect to $\Gamma_i$). Besides $S_{v_i}$ is a part of some
$v_i$-cycle, $S_u$ is a part of the $u$-cycle $S'$. Since the
elements $u, v_i$ of the group $A\ast B$ do not belong to conjugate
cyclic subgroups and the length of the path $S_u$ is greater than
$10s=10\max_{z\in\{\,u, v_1,\ldots, v_n\}}(l(z))$ the paths
$S_{v_i}$ and $S_u$ are different.

\begin{figure}
\centering
\includegraphics[15cm,12.5cm]{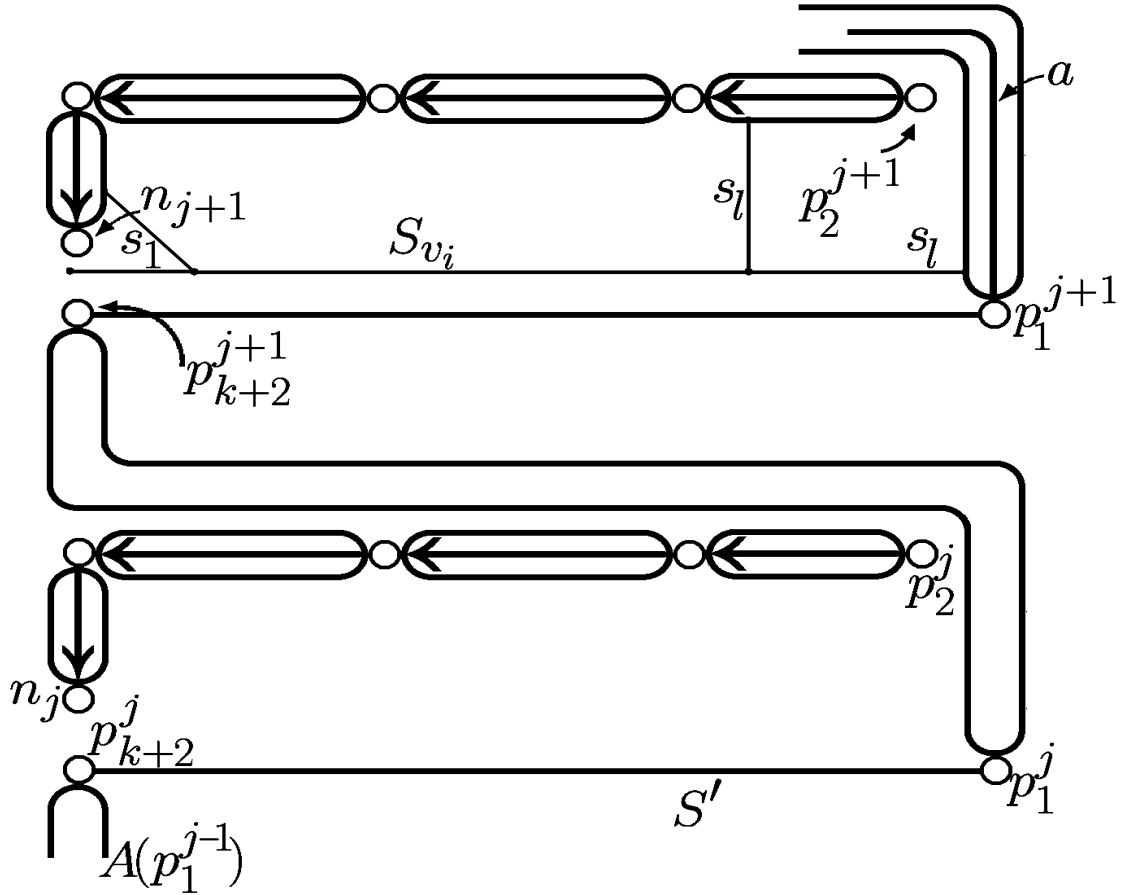}
\caption{The graph $\Delta$}
\end{figure}
\newpage

Suppose that the length of the path $S_{v_i}$ is less or equal than
$k+4=l(S_u)+3$. The paths $S_{v_i}$ and $S_u$ and perhaps several
edges whose number is less than 4 compose the loop. Let $g$ be the
label of this loop. Then $g$ is an element of the group $G$ whose
length is less or equal than $2l(S_u)+6=2k+8<10k$ and
$\varphi(g)=1$. But this contradicts the condition on $\varphi$ and
the set $Q$. Thus the length of the path $S_{v_i}$ is greater than
$k+4=l(S_u)+3$. By the symmetry we may also assume that the length
of the path $T\setminus S_{v_i}$ is greater than $k+4$: the
structure of the part of $T$ in $\Delta_{j+2}$ is the same as in
$\Delta_{j+1}$. But in this case $\rho(\alpha(S_{v_i}),
\omega(S_{v_i}))\leqslant\min(l(S_u)+3, l(S_{v_i}), l(T\setminus
S_{v_i}))=\min(k+4, l(S_{v_i}), l(T\setminus S_{v_i}))=k+4$, since
$l(S_{v_i}), l(T'\setminus R')>k+4$. So the $v_i$-cycle $T$
containing $S_{v_i}$ has $(k+4)$-near vertices. This also
contradicts the conditions on $\varphi$.

Thus it is proved that for each $i$ and for each $v_i$-cycle $T$ in
the graph $\Delta$ there exists $j, 1\leqslant j\leqslant m,$ such
that all vertices of $T$ are contained in $\Delta_j\cup\Delta_{j+1}$
(indices are modulo $m$). We deduce also that each $u$-cycle of
$\Delta$ which does not start at $p_1^1$ belongs to two subgraphs
$\Delta_{k_1}, \Delta_{k_1+1}$. This can be established by the same
way as it was shown that the analogical statement is true for
$v_i$-cycles.

Now we shall denote the obtained graph $\Delta$ for number $m$ as
$\Delta_m'$. Consider the set of graphs $\Lambda_m=\Delta_{3m}',
m=1, 2,..$.

We shall show now that $|\varphi_{\Lambda_m}(v_i)|$ equals some
constant number $K_i$ which does not depend on $m$. Let $R_{i, m}$
be the set of lengths of all $v_i$-cycles of $\Lambda_m$. The local
structure of $\Lambda_m$ is the same: using the above notations and
regarding that $\Lambda_m$ is the union of $\Delta_1,\ldots,
\Delta_{3m}$ it is obvious that the subgraphs
$\Delta_k\cup\Delta_{k+1}$ and $\Delta_l\cup\Delta_{l+1}$ are
isomorphic and do not depend on $m$. Hence $R_{i, m}$ coincides with
the set of lengths of $v_i$-cycles concentrated in
$\Delta_1\cup\Delta_2$ and thereby $R_{i, m}=R_{i, t}$ for all $m,
t$. The same reasonings are true for all $u$-cycles of $\Lambda_m$
except for the $u$-cycle whose length is the multiple of $3m$ so
$|\varphi_{\Lambda_m}(u)|=mK$ for some constant $K$ which does not
depend on $m$.

Proposition is proved and therefore the theorem is also proved.

\begin{center}
\large{Acknowledgements}
\end{center}

I am grateful to Anton A. Klyachko, Ashot Minasyan, Denis Osin and
Henry Wilton for valuable conversations and for information about
omnipotence.

\begin{center}
\large{References}
\end{center}

1. Klyachko, A. A. Equations over groups, quasivarieties, and a
residual property of a free group. \textsl{Journal of group theory}
\textbf{2}: 319-327, 1999.

2. Wise, Daniel T.  Subgroup separability of graphs of free groups
with cyclic edge groups. \textsl{Q. J. Math.} 51, No.1, 107-129
(2000). [ISSN 0033-5606; ISSN 1464-3847]

3. Jitendra Bajpai. Omnipotence of surface groups. Masters Thesis,
McGill University, 2007.

4. Henry Wilton. Virtual retractions, conjugacy separability and
omnipotence. J. Algebra 323 (2010), pp. 323-335.

5. Marshall Hall, Jr. Coset representations in free groups.
\textsl{Trans. Amer. Math. Soc.}, 67:421--432, 1948.

6. Yedynak, V. V. Separability with respect to order.
\textsl{Vestnik Mosk. Univ. Ser. I Mat. Mekh.} \textbf{3}: 56-58,
2006.

7. Yedynak, V. V. Multielement order separability in free products
of groups. \textsl{Communications in Algebra} 38(\textbf{3}): 3448
–- 3455, 2010.

\end{document}